\documentclass{amsart}

\usepackage{amsmath}
\usepackage{amsfonts}
\usepackage{amssymb}
\usepackage{amsxtra}

\newcommand{\defm}[1]{{\it #1}}
\newcommand{\R}{\mathbb{R}}
\newcommand{\subeq}[2]{\mathord{\underbrace{\mathop{#1}}_{#2}}}
\newcommand{\eqv}{\Leftrightarrow}
\newcommand{\impl}{\Rightarrow}

\newtheorem{theorem}{Theorem}
\newtheorem{corollary}{Corollary}
\newtheorem{remark}{Remark}

\begin{document}

\title{The ellipticity principle for selfsimilar polytropic potential 
flow}

\author{Volker Elling}
\address{Dr.-Maria-Moormann-Str. 12, D-48231 Warendorf, Germany}
\address{Dept. of Mathematics, Stanford University, 450 Serra Mall, Building 380\\Stanford, CA 94305-2125, USA}
\email{velling@stanford.edu}

\author{Tai-Ping Liu}

\address{Dept. of Mathematics, Stanford University, 450 Serra Mall, Building 380\\Stanford, CA 94305-2125, USA}
\email{liu@math.stanford.edu}

\begin{abstract}
    We consider self-similar potential flow for compressible gas with polytropic pressure law. 
    Self-similar solutions arise as large-time asymptotes of general solutions, and as 
    exact solutions of many important special cases like Mach reflection, 
    multidimensional Riemann problems, or flow around corners.

    Self-similar potential flow is a quasilinear second-order PDE of mixed type which is hyperbolic at infinity 
    (if the velocity is globally bounded).
    The type in each point is determined by the local pseudo-Mach number $L$, with $L<1$ resp.\ $L>1$ corresponding 
    to elliptic resp.\ hyperbolic regions. We prove an \defm{ellipticity principle}: the interior of 
    a parabolic-elliptic region of a sufficiently smooth solution must be elliptic; in fact $L$ must 
    be bounded above away from $1$ by a domain-dependent function. In particular there are no open 
    parabolic regions. We also discuss the case of slip boundary conditions at straight solid walls.
\end{abstract}

\maketitle


\section{Introduction}

We consider the isentropic Euler equations of compressible gas dynamics in $d$ space dimensions:
\begin{alignat}{1}
    \rho_t + \nabla\cdot(\rho\vec v) &= 0 \label{eq:rho} \\
    (\rho\vec v)_t + \nabla\cdot(\rho\vec v\otimes\vec v) + \nabla(p(\rho)) &= 0. \label{eq:mom}
\end{alignat}
Hereafter, $\nabla$ denotes the gradient with respect either to the space coordinates ${\vec x}=(x^1,x^2,\cdots,x^d)$ or
to similarity coordinates $t^{-1}\vec x$.
Throughout this paper, superscripts represent vector components whereas subscripts represent partial derivatives.
${\vec v}=(v^1,v^2,\cdots,v^d)$ is the velocity of the gas, $\rho$ the density, and $p=p(\rho)$ follows the polytropic pressure law
\begin{alignat}{1}
    p(\rho) &= \frac{c_0^2\rho_0}{\gamma}\left(\frac{\rho}{\rho_0}\right)^\gamma \label{eq:p-polytropic}
\end{alignat}
for $\gamma\neq 0$ (here $c_0$ is the sound speed at density $\rho_0$).

For smooth solutions, substituting \eqref{eq:rho} into \eqref{eq:mom} yields the simpler form
\begin{alignat}{1}
    \vec v_t+\nabla\vec v\cdot\vec v + \nabla(\pi(\rho)) &= 0. \tag{\ref{eq:mom}'}\label{eq:v}
\end{alignat}
Here $\pi$ is defined as
\begin{alignat}{1}
    \pi(\rho) &= \begin{cases}
	\frac{c_0^2}{\gamma-1}(\rho/\rho_0)^{\gamma-1} = \frac{1}{\gamma-1}c^2(\rho), & \gamma\neq 1 \\
	c_0^2\log\frac{\rho}{\rho_0}, & \gamma=1.
    \end{cases}
\end{alignat}

The Euler equations possess highly unstable vortex sheets. To focus on the acoustic waves and shock waves, one often considers the \defm{irrotationality} assumption $v^i_j=v^j_i$ (where $i,j=1,\dotsc,d$), which leads to 
\begin{alignat*}{1}
    \vec v &= \nabla_{\vec x}\phi
\end{alignat*}
for some scalar \defm{potential} function $\phi$. For smooth flows, substituting this into \eqref{eq:v} yields, for $i=1,\dotsc,d$,
\begin{alignat*}{1}
    0 
    &= \phi_{it} + \nabla\phi_i\cdot\nabla\phi + \pi(\rho)_i \\
    &= (\phi_t + \frac{|\nabla\phi|^2}{2} + \pi(\rho))_i.
\end{alignat*}
Thus, for some constant $A$, 
\begin{alignat*}{1}
    \rho &= \pi^{-1}(A-\phi_t-\frac{|\nabla\phi|^2}{2}).
\end{alignat*}
Substituting this into \eqref{eq:rho} yields a single second-order quasilinear hyperbolic equation, the \defm{potential flow} equation, for a scalar field $\phi$:
using $c^2=p_\rho$ and $(\pi^{-1})'=\frac{\rho}{c^2}$, we get
\begin{alignat}{1}
    \phi_{tt} + 2\nabla\phi_t\cdot\nabla\phi + \sum_{i,j=1}^d\phi_i\phi_j\phi_{ij} - c^2\Delta\phi &= 0. \label{eq:potential-flow}
\end{alignat}

Due to the invariance of first-order systems of conservation laws under the scaling $x\leftarrow sx$, 
$t\leftarrow st$ for $s>0$, there is an important class of solutions to \eqref{eq:rho} and \eqref{eq:v}, 
the \defm{self-similar solutions}, which depend only on the similarity coordinates
$$\vec\xi:=\frac{\vec x}{t}.$$
Self-similar solutions arise as large-time asymptotes of unsteady solutions and as exact solutions 
to many important special cases like Riemann problems, Mach reflection, flow around solid corners etc. 

For potential flow the appropriate self-similar ansatz is
\begin{alignat*}{1}
    \phi(t,\vec x) &= t\psi(\vec\xi)
\end{alignat*}
because $\vec v = \nabla_{\vec x}\phi = \nabla_{\vec\xi}\psi$.
\eqref{eq:potential-flow} turns into the self-similar potential flow equation:
\begin{alignat}{1}
    c^2\Delta\psi-\sum_{i,j=1}^d(\psi_i-\xi^i)(\psi_j-\xi^j)\psi_{ij} &= 0. \label{eq:psi}
\end{alignat}
The convenient change of variables
\begin{alignat*}{1}
    \chi &:= \psi-\frac{|\vec\xi|^2}{2}
\end{alignat*}
yields an equivalent equation whose coefficients do not directly depend on $\vec\xi$:
\begin{alignat}{1}
    c^2\Delta\chi-\sum_{i,j=1}^d\chi_i\chi_j\chi_{ij} &= |\nabla\chi|^2-dc^2 \label{eq:chi}
\end{alignat}

For polytropic pressure laws, the sound speed $c$ has a particularly simple form. For $\gamma=1$ (isothermal flow),
$c=c_0$; for $\gamma\neq 1$,
\begin{alignat}{1}
    c^2 &= (\gamma-1)(A-\chi-\frac{|\nabla\chi|^2}{2}), \label{eq:css}
\end{alignat}
Clearly, adding a constant to $\psi$ resp.\ $\chi$ does not matter in 
\eqref{eq:psi} resp.\ \eqref{eq:chi}, so we may take $A=0$ in \eqref{eq:css} for simplicity:
\begin{alignat}{1}
    c^2 &= (1-\gamma)(\chi+\frac{|\nabla\chi|^2}{2}). \label{eq:cs}
\end{alignat}

$\nabla\chi=\nabla\psi-\vec\xi$ is called \defm{pseudo-velocity}. 

\begin{remark}
	\label{rem:symmetries}%
	\eqref{eq:chi} inherits a number of symmetries from \eqref{eq:rho}, \eqref{eq:mom}: 
	\begin{enumerate}
	\item it is invariant under translation (which corresponds to the transformations $v\leftarrow v+v_0$, $x\leftarrow x-v_0t$ 	(with parameter $v_0\in\R^d$) in $(t,x)$ coordinates),
	\item it is invariant under rotation, and
	\item it is invariant under velocity scaling (i.e.\ if $\chi$ is a global solution, 
	then $\tilde\chi(\vec\xi):=s^2\chi(s^{-1}\vec\xi)$ is a global solution as well, for any $s>0$).
	\end{enumerate}
\end{remark}

For many interesting flow patterns, \eqref{eq:chi} is of mixed type. The type
is determined by the (local) \defm{pseudo-Mach number}
\begin{alignat*}{1}
    L &:= \frac{|\nabla\chi|}{c},
\end{alignat*}
with $0\leq L<1$ for elliptic (pseudo-subsonic), $L=1$ for parabolic (pseudo-sonic), $L>1$ for hyperbolic (pseudo-supersonic) regions.

It is important to identify the elliptic and hyperbolic regions because their mathematical and physical 
properties are quite different. The equations \eqref{eq:chi} and \eqref{eq:cs} are highly nonlinear, so the 
regions are not known a priori from the flow data. 
In contrast to steady potential flow which has many interesting solutions of a single type, any 
self-similar potential flow solution in the entire space with bounded velocity must be hyperbolic at infinity
and most cases appear to have an elliptic region as well.

Our main result, Theorem \ref{th:ellipticity}, states that if the local pseudo-Mach 
number of a $C^3$ solution of \eqref{eq:chi} is
$\leq 1$ on the smooth boundary of a bounded domain $\Omega$, then 
it is bounded above away from $1$ in the interior by a barrier function which depends 
only on the domain and on an upper bound on the sound speed. 
Such an example is useful for a variety of reasons; in the construction of self-similar 
solutions it guarantees that perturbing the smooth boundary of a 
bounded elliptic region in a way that keeps $L\leq 1$ on the boundary \emph{must}\footnote{It is 
expected that for boundary conditions that admit a unique 
solution, the solution is $C^3$ (in fact analytic) in the interior and that continuous 
perturbations of the boundary and boundary data yield continuous changes of the solution in any $C^k$ norm.} keep the 
solution elliptic inside. In absence of such a result, 
a possibly infinite number of parabolic-hyperbolic ``bubbles'' could arise in 
the interior which would make analysis prohibitively complicated.

A variety of maximum principles for fluid variables has appeared in the 
literature. A direct analogue is the ellipticity principle derived in \cite{yuxi-zheng-ellipticity} (see also 
\cite{li-zhang-yang,yuxi-zheng-book}) for the self-similar pressure-gradient equations.
\cite{serre-p-maxprinciple} shows that, under certain conditions, 
pressure in steady inviscid incompressible flow has extrema only at rest points. 

In \cite[Chapter 15]{gilbarg-trudinger}, rather general maximum principles for the gradient of second-order quasilinear elliptic PDE are provided. They apply in particular to equations of the form
$$a^{ij}(\nabla u)\partial_{ij}u+b^i(\nabla u)\partial_iu=0,$$
i.e.\ in absence of $0$th order terms and with coefficients depending on the gradient only. Among the corollaries
are strong maximum principles for the Mach number (or equivalently density or velocity) in \emph{steady} potential flow
as well as for the pseudo-Mach number $L$ in the isothermal case ($\gamma=1$) for self-similar potential flow. 
However, the case $\gamma\neq 1$ is not covered.

\section{The ellipticity principle for self-similar polytropic potential flows}

\begin{theorem}[Ellipticity principle]
    \label{th:ellipticity}%
    Let $d\geq 2$ and $\gamma>-1$.
    Let $\Omega\subset\R^d$ open and bounded. There is a positive constant $\delta$ so that, 
    for any $\hat c>0$ and any $b\in C^2(\Omega)$ with $|\nabla b|\leq\frac{\delta}{\hat c}$,                 $|\nabla^2b|\leq\frac{\delta}{\hat c^2}$ and for any solution 
    $\chi\in C^3(\Omega)$ of \eqref{eq:chi} with $L\leq 1$, $\rho>0$ and $c\leq\hat c$, 
    either 
    $$L^2\leq 1-\delta$$
    or
    $L^2+b$ does not attain its maximum in $\Omega$.
\end{theorem}
\begin{remark}
    At its heart the bound depends only on the domain, since the velocities $c,\vec v,\vec\xi$ 
    can be rescaled (see Remark \ref{rem:symmetries}) to allow $\hat c=1$.
\end{remark}
\begin{proof}
    During the proof, all equations and inequalities are meant to hold
    in the maximum point only. 
    The $O$ notation is with respect to $\delta\downarrow 0$ and uniform in $c$.

    Assume that $L^2+b$ has a maximum in some interior point with $L^2>1-\delta$.
    We fix $\delta\leq 1$ right away, so we may assume $|\nabla\chi|=L c\neq 0$.
    Since the equations are rotation-invariant, we may, without loss of generality, choose $\chi_1>0$ and
    $\chi_j=0$ for all $j>1$. 
    The first-order conditions for a maximum imply
    \begin{alignat*}{1}
        0 = (L^2+b)_1 &= \frac{2\chi_1\chi_{11}c^2+(\gamma-1)\chi_1^3(1+\chi_{11})}{c^4}+b_1 \\
        \eqv \chi_{11} &= \frac{-cb_1-(\gamma-1)L^3}{L(2+(\gamma-1)L^2)} = \frac{1-\gamma}{\gamma+1}+O(\delta).
    \end{alignat*}
    For $j>1$,
    \begin{alignat*}{1}
        0 = (L^2+b)_j &= \frac{2\chi_1\chi_{1j}c^2+(\gamma-1)\chi_1^3\chi_{1j}}{c^4}+b_2 \\
        \eqv \chi_{1j} &= O(\delta).
    \end{alignat*}
    From \eqref{eq:chi},
    \begin{alignat}{1}
        \sum_{j>1}\chi_{jj} &= (L^2-1)(\chi_{11}+1)+L^2-d = 1-d + O(\delta) \label{eq:chi_jj}
    \end{alignat}
    Moreover, using $(f/g)''=(f''g^2-fgg''+2f(g')^2-2f'gg')/g^3$, the second-order conditions yield
    \begin{alignat*}{1}
        0 &\geq (L^2+b)_{11} \\
	&= \frac{2(\chi_1\chi_{111}+\chi_{11}^2+\sum_{j>1}\chi_{1j}^2)c^4}{c^6} \\
	&+\frac{(\gamma-1)\chi_1^2c^2(\chi_{11}(1+\chi_{11})+\sum_{j>1}\chi_{1j}^2+\chi_1\chi_{111})}{c^6} \\
        &+\frac{2(\gamma-1)^2\chi_1^4(1+\chi_{11})^2}{c^6}+\frac{4(\gamma-1)\chi_1^2c^2\chi_{11}(1+\chi_{11})}{c^6} + c^{-2}O(\delta)\\
        &= c^{-1}L(2+(\gamma-1)L^2)\chi_{111}\\
	&+c^{-2}(\chi_{11}+2(\gamma-1)L^2(1+\chi_{11}))(2\chi_{11}+(\gamma-1)L^2(1+\chi_{11}))\\
        &+c^{-2}(2+(\gamma-1)L^2)\sum_{j>1}\chi_{1j}^2 + c^{-2}O(\delta) \\
        \eqv\chi_{111} &\leq \frac{1}{cL(2+(\gamma-1)L^2)}
        \Big(\subeq{(\chi_{11}+2(\gamma-1)L^2(1+\chi_{11}))}{=O(1)}
        \subeq{(2\chi_{11}+(\gamma-1)L^2(1+\chi_{11}))}{=O(\delta)} \\
        &+\subeq{(2+(\gamma-1)L^2)}{=O(1)}\subeq{\sum_{j>1}\chi_{1j}^2}{=O(\delta)} 
        +O(\delta)\Big) = c^{-1}O(\delta) 
    \end{alignat*}
    For $j>1$,
    \begin{alignat}{1}
        0 &\geq (L^2+b)_{jj} \notag\\
	&= \frac{2(\chi_1\chi_{1jj}+\sum_{i=1}^d\chi_{ij}^2)c^4}{c^6}
	+\frac{(\gamma-1)\chi_1^2c^2(\chi_{jj}+\sum_{i=1}^d\chi_{ij}^2+\chi_1\chi_{1jj})}{c^6}\notag\\
        &+\frac{2(\gamma-1)^2\chi_1^4\chi_{1j}^2}{c^6}
        +\frac{4(\gamma-1)\chi_1^2\chi_{1j}^2c^2}{c^6} + c^{-2}O(\delta) \notag \\
        &= \frac{L(2+(\gamma-1)L^2)}{c}\chi_{1jj}+\frac{(2+(\gamma-1)L^2)(1+2(\gamma-1)L^2)}{c^2}\chi_{1j}^2\notag\\
	&+\frac{2+(\gamma-1)L^2}{c^2}\sum_{k>1}\chi_{kj}^2 + \frac{(\gamma-1)L^2}{c^2}\chi_{jj} + c^{-2}O(\delta) \notag 
    \end{alignat}
    \begin{alignat*}{1}
        \impl\quad c\chi_{1jj} \leq \subeq{\frac{-1}{L}}{=-1+O(\delta)}\Big(\subeq{(1+2(\gamma-1)L^2)\chi_{1j}^2}{\geq 0}
         +\subeq{\sum_{k>1}\chi_{kj}^2}{\geq\chi_{jj}^2}\Big)
          -\subeq{\frac{(\gamma-1)L}{2+(\gamma-1)L^2}}{=(\gamma-1)/(\gamma+1)+O(\delta)}\chi_{jj}+O(\delta) 
    \end{alignat*}
    \begin{alignat}{1}
	\impl\qquad c\sum_{j>1}\chi_{1jj} &\leq -\sum_{j>1}\chi_{jj}^2 \notag\\
	&+\frac{\gamma-1}{\gamma+1}(d-1)+O(\delta) \label{eq:chi122-1a}
   \end{alignat}
    Minimizing $\sum_{j>1}\chi_{jj}^2$ with respect to (see \eqref{eq:chi_jj})
    $$\sum_{j>1}\chi_{jj} = 1-d+O(\delta);$$
    the minimum is 
    $$\chi_{22}=\dotsc=\chi_{dd}=\frac{1-d+O(\delta)}{d-1}=-1+O(\delta);$$
    substituting this into \eqref{eq:chi122-1a} yields
    \begin{alignat}{1}
        c\sum_{j>1}\chi_{1jj} &\leq \frac{-2(d-1)}{\gamma+1}+O(\delta) \label{eq:chi122-1}
    \end{alignat}
    At this point it is necessary to take $\partial_1$ of \eqref{eq:chi}:
    \begin{alignat*}{1}
        & \subeq{(c^2-\chi_1^2)}{\geq 0}\subeq{\chi_{111}}{\leq c^{-1}O(\delta)}+\sum_{j>1}(c^2-\subeq{\chi_j^2}{=0})\chi_{1jj}-2\sum_{j\neq k}\subeq{\chi_j\chi_k}{=0}\chi_{1jk} \\
        & + \subeq{((1-\gamma)\chi_1+(1-\gamma)(\chi_1\subeq{\chi_{11}}{=\frac{1-\gamma}{\gamma+1}+O(\delta)}+\sum_{j>1}\subeq{\chi_j}{=0}\chi_{1j})-2\chi_1\subeq{\chi_{11}}{=\frac{1-\gamma}{\gamma+1}+O(\delta)})}{=O(\delta)\chi_1}\subeq{\chi_{11}}{=O(\delta)\chi_1} \\
        & + \sum_{k>1}((1-\gamma)\chi_1+(1-\gamma)(\chi_1\subeq{\chi_{11}}{=\frac{1-\gamma}{\gamma+1}+O(\delta)}+\sum_{j>1}\subeq{\chi_j}{=0}\chi_{1j})-2\subeq{\chi_k}{=0}\chi_{1k})\chi_{kk} \\
        & -2\sum_{j>1}(\chi_{11}\chi_j+\chi_1\chi_{1j})\subeq{\chi_{1j}}{=O(\delta)} 
        -2\sum_{1<j<k}(\chi_{1j}\subeq{\chi_k}{=0}+\subeq{\chi_j}{=0}\chi_{1k})\chi_{jk} \\
        &=d(\gamma-1)\chi_1+(2+d(\gamma-1))(\chi_1\subeq{\chi_{11}}{=\frac{1-\gamma}{\gamma+1}+O(\delta)}+\sum_{k>1}\subeq{\chi_k}{=0}\chi_{1k}) 
    \end{alignat*}
    \begin{alignat*}{1}
        \Leftrightarrow\quad cO(\delta)+2\frac{\gamma-1}{\gamma+1}\chi_1(\sum_{j>1}\chi_{jj}+d-1) &\leq c^2\sum_{j>1}\chi_{1jj} 
    \end{alignat*}
    Using \eqref{eq:chi_jj} we see
    \begin{alignat}{1}
        c\sum_{j>1}\chi_{1jj} &\geq O(\delta) \label{eq:chi122-2}
    \end{alignat}
    Comparing \eqref{eq:chi122-1} and \eqref{eq:chi122-2} yields a contradiction (for $d\geq 2$) 
    for $\delta$ chosen sufficiently small.
\end{proof}

\begin{corollary}
    \label{th:open-sonic-regions}%
    For $d>1$, \eqref{eq:chi} does not admit solutions $\chi\in C^3$ that are parabolic on some open set.
\end{corollary}

\begin{remark}
    Corollary \ref{th:open-sonic-regions} does not hold for $d=1$. 
    In this case \eqref{eq:chi} reduces to the second-order ODE 
    \begin{alignat}{1}
        (c^2-(\chi')^2)\chi'' &= (\chi')^2-c^2,\qquad c^2=\begin{cases}(1-\gamma)(\chi+\frac{1}{2}(\chi')^2), & \gamma\neq 1\\ \text{const}, & \gamma=1\end{cases}.\label{eq:chi-1d}
    \end{alignat}
    Solutions of the first-order ODE $c^2=(\chi')^2$ 
    satisfies \eqref{eq:chi-1d} as well, but they are parabolic everywhere.
    (Note that those solutions correspond to rarefaction waves.) 

    Nevertheless, solutions to \eqref{eq:chi-1d} are affine on any 
    non-parabolic interval (this is obvious from \eqref{eq:psi} for $d=1$),
    so a weaker version of ellipticity holds trivially.
\end{remark}

\section{Straight solid walls with slip condition}

Theorem \ref{th:ellipticity} by itself is suitable for parabolic boundaries (boundary condition $L=1$)
or regularizations ($L=1-\epsilon$). Here we extend it to another common boundary condition, the slip condition
\begin{alignat}{1}
	\frac{\partial\chi}{\partial n} &= 0 \label{eq:slip}
\end{alignat}
at a solid wall. We consider straight walls only because in self-similar flow curved walls are not physically interesting.

\begin{theorem}
    \label{th:wall}%
    Let $d\geq 2$ and $\gamma>-1$.
    Let $U$ be a small open ball centered on $\{0\}\times\R^{d-1}$ 
    and set $\Omega:=U\cap((0,\infty)\times\R^{d-1})$ and $\Gamma:=U\cap(\{0\}\times\R^{d-1})$.
    There is a positive constant $\delta$ so that, for any $\hat c>0$, any $b\in C^2(\Omega\cup\Gamma)$
    with $\frac{\partial b}{\partial n}=0$ on $\Gamma$,
    $|\nabla b|\leq\frac{\delta}{\hat c}$, $|\nabla^2b|\leq\frac{\delta}{\hat c^2}$ in $\Omega$, 
    and any $\chi\in C^3(\Omega\cup\Gamma)$ 
    that solves \eqref{eq:chi} in $\Omega$, satisfies \eqref{eq:slip} on $\Gamma$
    and has $L\leq 1$, $\rho>0$ and $c\leq\hat c$ in $\Omega\cup\Gamma$, either
    $$L^2\leq 1-\delta$$
    or $L^2+b$ does not attain its maximum in $\Omega\cup\Gamma$.
\end{theorem}
\begin{proof}
    In this setting \eqref{eq:slip} means $\chi_1=0$ on $\Gamma$. Taking tangential derivatives in 
    \eqref{eq:slip} yields
    $$\chi_{1i}=0,\quad\chi_{1ij}=0\qquad\forall i,j\in\{2,\dotsc,d\}$$
    on $\Gamma$.
    This implies
    $$(c^2)_1=(1-\gamma)(\subeq{\chi_1}{=0}+\subeq{\chi_1}{=0}\chi_{11}+\sum_{j=2}^d\chi_j\subeq{\chi_{1j}}{=0})=0$$
    on $\Gamma$ as well.
    Moreover take $\partial_1$ of \eqref{eq:chi}:
    \begin{alignat*}{1}
	&c^2\chi_{111}
        +\subeq{(c^2)_1}{=0}\chi_{11} 
        +\sum_{i=2}^d(c^2\subeq{\chi_{1ii}}{=0}
                     +\subeq{(c^2)_1}{=0}\chi_{ii}) 
        -\subeq{\chi_1^2}{=0}\chi_{111}-2\subeq{\chi_1}{=0}\chi_{11}^2 \\
        &-2\sum_{j=2}^d(\subeq{\chi_1}{=0}\chi_j\chi_{11j}+(\chi_1\chi_j)_1\subeq{\chi_{1j}}{=0})
        -\sum_{i,j=2}^d(\subeq{\chi_{1i}}{=0}\chi_j+\chi_i\subeq{\chi_{1j}}{=0}+\chi_i\chi_j\subeq{\chi_{1ij}}{=0}) \\
        &= 2\subeq{\chi_1}{=0}\chi_{11}+2\sum_{j=2}^d\chi_j\subeq{\chi_{1j}}{=0}-d\subeq{(c^2)_1}{=0} \qquad\text{on $\Gamma$}
    \end{alignat*}
    \begin{alignat*}{1}
	\Rightarrow \chi_{111} &= 0 \qquad\text{on $\Gamma$}
    \end{alignat*}
    Combining all results we see that even reflection
    $$\chi(\xi^1,\xi^2,\dotsc,\xi^d):=\chi(-\xi^1,\xi^2,\dotsc,\xi^d)$$
    yields a $C^3(U)$ extension of $\chi$ that satisfies \eqref{eq:chi} in all of $U$. 
    Moreover since $b_1=0$ on $\Gamma$, we can extend $b$ in the same manner to $C^2(U)$. All conditions of Theorem
    \ref{th:ellipticity} are satisfied, so it implies the desired result.
\end{proof}

\begin{remark}
	Any solid wall can be rotated and translated (see Remark \ref{rem:symmetries})
        to coincide with $\Gamma$ in 
        Theorem \ref{th:wall}.
\end{remark}

\section*{Acknowledgements}

This material is based upon work supported by the National Science Foundation under Grant no. DMS 0104019. 
The first author was supported by an SAP/Stanford Graduate Fellowship. The authors would like to thank Marshall Slemrod for 
valuable discussions and for his interest in this work.

\bibliographystyle{amsalpha}
\bibliography{../thesis/elling}

\end{document}